\documentclass[12pt]{article}

\usepackage{amsmath,amssymb,amsthm,amscd,a4wide}
\usepackage{hyperref}
\hypersetup{colorlinks,citecolor=blue,filecolor=black,linkcolor=blue,urlcolor=blue}

\begin{document}

\newcommand{\End}{{\rm{End}\ts}}
\newcommand{\Hom}{{\rm{Hom}}}
\newcommand{\Mat}{{\rm{Mat}}}
\newcommand{\ad}{{\rm{ad}\ts}}
\newcommand{\ch}{{\rm{ch}\ts}}
\newcommand{\chara}{{\rm{char}\ts}}
\newcommand{\diag}{ {\rm diag}}
\newcommand{\pr}{^{\tss\prime}}
\newcommand{\non}{\nonumber}
\newcommand{\wt}{\widetilde}
\newcommand{\wh}{\widehat}
\newcommand{\ot}{\otimes}
\newcommand{\la}{\lambda}
\newcommand{\ls}{\ts\lambda\ts}
\newcommand{\La}{\Lambda}
\newcommand{\De}{\Delta}
\newcommand{\al}{\alpha}
\newcommand{\be}{\beta}
\newcommand{\ga}{\gamma}
\newcommand{\Ga}{\Gamma}
\newcommand{\ep}{\epsilon}
\newcommand{\ka}{\kappa}
\newcommand{\vk}{\varkappa}
\newcommand{\vt}{\vartheta}
\newcommand{\si}{\sigma}
\newcommand{\vp}{\varphi}
\newcommand{\de}{\delta}
\newcommand{\ze}{\zeta}
\newcommand{\om}{\omega}
\newcommand{\ee}{\epsilon^{}}
\newcommand{\su}{s^{}}
\newcommand{\hra}{\hookrightarrow}
\newcommand{\ve}{\varepsilon}
\newcommand{\ts}{\,}
\newcommand{\vac}{\mathbf{1}}
\newcommand{\vacf}{\tss|0\rangle}
\newcommand{\di}{\partial}
\newcommand{\qin}{q^{-1}}
\newcommand{\tss}{\hspace{1pt}}
\newcommand{\Sr}{ {\rm S}}
\newcommand{\U}{ {\rm U}}
\newcommand{\BL}{ {\overline L}}
\newcommand{\BE}{ {\overline E}}
\newcommand{\BP}{ {\overline P}}
\newcommand{\AAb}{\mathbb{A}\tss}
\newcommand{\CC}{\mathbb{C}\tss}
\newcommand{\KK}{\mathbb{K}\tss}
\newcommand{\QQ}{\mathbb{Q}\tss}
\newcommand{\SSb}{\mathbb{S}\tss}
\newcommand{\ZZ}{\mathbb{Z}\tss}
\newcommand{\X}{ {\rm X}}
\newcommand{\Y}{ {\rm Y}}
\newcommand{\Z}{{\rm Z}}
\newcommand{\Ac}{\mathcal{A}}
\newcommand{\Lc}{\mathcal{L}}
\newcommand{\Mc}{\mathcal{M}}
\newcommand{\Pc}{\mathcal{P}}
\newcommand{\Qc}{\mathcal{Q}}
\newcommand{\Tc}{\mathcal{T}}
\newcommand{\Sc}{\mathcal{S}}
\newcommand{\Bc}{\mathcal{B}}
\newcommand{\Ec}{\mathcal{E}}
\newcommand{\Fc}{\mathcal{F}}
\newcommand{\Hc}{\mathcal{H}}
\newcommand{\Uc}{\mathcal{U}}
\newcommand{\Vc}{\mathcal{V}}
\newcommand{\Wc}{\mathcal{W}}
\newcommand{\Yc}{\mathcal{Y}}
\newcommand{\Ar}{{\rm A}}
\newcommand{\Br}{{\rm B}}
\newcommand{\Ir}{{\rm I}}
\newcommand{\Fr}{{\rm F}}
\newcommand{\Jr}{{\rm J}}
\newcommand{\Or}{{\rm O}}
\newcommand{\Tr}{{\rm T}}
\newcommand{\GL}{{\rm GL}}
\newcommand{\Spr}{{\rm Sp}}
\newcommand{\Rr}{{\rm R}}
\newcommand{\Zr}{{\rm Z}}
\newcommand{\gl}{\mathfrak{gl}}
\newcommand{\middd}{{\rm mid}}
\newcommand{\ev}{{\rm ev}}
\newcommand{\Pf}{{\rm Pf}}
\newcommand{\Norm}{{\rm Norm\tss}}
\newcommand{\oa}{\mathfrak{o}}
\newcommand{\spa}{\mathfrak{sp}}
\newcommand{\osp}{\mathfrak{osp}}
\newcommand{\g}{\mathfrak{g}}
\newcommand{\h}{\mathfrak h}
\newcommand{\n}{\mathfrak n}
\newcommand{\z}{\mathfrak{z}}
\newcommand{\Zgot}{\mathfrak{Z}}
\newcommand{\p}{\mathfrak{p}}
\newcommand{\sll}{\mathfrak{sl}}
\newcommand{\agot}{\mathfrak{a}}
\newcommand{\qdet}{ {\rm qdet}\ts}
\newcommand{\Ber}{ {\rm Ber}\ts}
\newcommand{\HC}{ {\mathcal HC}}
\newcommand{\cdet}{ {\rm cdet}}
\newcommand{\tr}{ {\rm tr}}
\newcommand{\gr}{ {\rm gr}}
\newcommand{\str}{ {\rm str}}
\newcommand{\loc}{{\rm loc}}
\newcommand{\Gr}{{\rm G}}
\newcommand{\sgn}{ {\rm sgn}\ts}
\newcommand{\ba}{\bar{a}}
\newcommand{\bb}{\bar{b}}
\newcommand{\bi}{\bar{\imath}}
\newcommand{\bj}{\bar{\jmath}}
\newcommand{\bk}{\bar{k}}
\newcommand{\bl}{\bar{l}}
\newcommand{\hb}{\mathbf{h}}
\newcommand{\Sym}{\mathfrak S}
\newcommand{\fand}{\quad\text{and}\quad}
\newcommand{\Fand}{\qquad\text{and}\qquad}
\newcommand{\For}{\qquad\text{or}\qquad}
\newcommand{\OR}{\qquad\text{or}\qquad}

\newcommand{\germ}{\mathfrak}
\newcommand{\cprime}{$'$}
\newcommand{\on}{\operatorname}
\newcommand{\+}{\mathop{\oplus}}
\newcommand{\mc}{\mathcal}
\newcommand{\mf}{\mathfrak}
\newcommand{\mb}{\mathbb}
\newcommand{\fing}{\mf{g}}
\newcommand{\finq}{\mf{q}}
\newcommand{\finl}{\mf{l}}
\newcommand{\finb}{\mf{b}}
\newcommand{\finm}{\mathfrak{m}}
\newcommand{\affg}{\widehat{\mf{gl}}_N}
\newcommand{\affq}{\widehat{\mf{q}}}
\newcommand{\affb}{\widehat{\mf{b}}}
\newcommand{\affa}{\widehat{\mf{a}}}
\newcommand{\isomap}{{\;\stackrel{_\sim}{\to}\;}}
\newcommand{\N}{\mathbb{N}}
\newcommand{\Q}{\mathbb{Q}}
\newcommand{\ra}{\rightarrow}
\newcommand{\lam}{\lambda}

\renewcommand{\theequation}{\arabic{section}.\arabic{equation}}

\newtheorem{thm}{Theorem}[section]
\newtheorem{lem}[thm]{Lemma}
\newtheorem{prop}[thm]{Proposition}
\newtheorem{cor}[thm]{Corollary}
\newtheorem{conj}[thm]{Conjecture}
\newtheorem*{mthm}{Main Theorem}
\newtheorem*{mthma}{Theorem A}
\newtheorem*{mthmb}{Theorem B}

\theoremstyle{definition}
\newtheorem{defin}[thm]{Definition}

\theoremstyle{remark}
\newtheorem{remark}[thm]{Remark}
\newtheorem{example}[thm]{Example}

\newcommand{\bth}{\begin{thm}}
\renewcommand{\eth}{\end{thm}}
\newcommand{\bpr}{\begin{prop}}
\newcommand{\epr}{\end{prop}}
\newcommand{\ble}{\begin{lem}}
\newcommand{\ele}{\end{lem}}
\newcommand{\bco}{\begin{cor}}
\newcommand{\eco}{\end{cor}}
\newcommand{\bde}{\begin{defin}}
\newcommand{\ede}{\end{defin}}
\newcommand{\bex}{\begin{example}}
\newcommand{\eex}{\end{example}}
\newcommand{\bre}{\begin{remark}}
\newcommand{\ere}{\end{remark}}
\newcommand{\bcj}{\begin{conj}}
\newcommand{\ecj}{\end{conj}}

\newcommand{\bal}{\begin{aligned}}
\newcommand{\eal}{\end{aligned}}
\newcommand{\beq}{\begin{equation}}
\newcommand{\eeq}{\end{equation}}
\newcommand{\ben}{\begin{equation*}}
\newcommand{\een}{\end{equation*}}

\newcommand{\bpf}{\begin{proof}}
\newcommand{\epf}{\end{proof}}

\def\beql#1{\begin{equation}\label{#1}}

\title{\Large\bf Explicit generators in rectangular affine\\ $\Wc$-algebras
of type $A$}

\author{{Tomoyuki Arakawa \quad and\quad Alexander Molev}}

\date{} 
\maketitle

\vspace{25 mm}

\begin{abstract}
We produce in an explicit form free generators of the affine $\Wc$-algebra
of type $A$ associated with a nilpotent matrix whose
Jordan blocks are of the same size. This includes the principal nilpotent case
and we thus recover the quantum Miura transformation of Fateev and Lukyanov.

\vspace{5 mm}

{\it Mathematics Subject Classification 2010\tss}: 17B69

\vspace{5 mm}

{\it Key words\tss}: vertex algebra, affine W-algebra, Miura transformation

\end{abstract}

\vspace{45 mm}

\noindent
Research Institute for Mathematical Sciences\newline
Kyoto University, Kyoto 606-8502, Japan\newline
arakawa@kurims.kyoto-u.ac.jp

\vspace{7 mm}

\noindent
School of Mathematics and Statistics\newline
University of Sydney, NSW 2006, Australia\newline
alexander.molev@sydney.edu.au

\newpage

\section{Introduction}
\label{sec:int}
\setcounter{equation}{0}

Let $\g$ be a reductive Lie algebra over $\CC$
equipped with a symmetric invariant bilinear form
and let $f$ be a nilpotent element of $\g$.
The corresponding {\em affine
$\Wc$-algebra $\Wc^{\tss k}(\g,f)$ at the level $k\in\CC$} is defined by
the generalized quantized Drinfeld--Sokolov reduction; see
\cite{ff:qd}, \cite{krw:qr} and \cite{kw:qr}. The case of a principal
nilpotent element $f$ was known much earlier and the corresponding $\Wc$-algebras
were intensively studied; see \cite[Ch.~15]{fb:va} for a detailed review of their structure,
and \cite{a:rtw} and \cite{a:rw} for representation theory.

In this paper we take $\g=\gl_N$. The Jordan type of a nilpotent element
$f\in\gl_N$ is a partition of $N$. We will work with the elements $f$
corresponding to partitions of the form $(l^{\tss n})$ so that the associated
Young diagram is the $n\times l$ rectangle with $n\tss l=N$.  Our main result
is an explicit construction of free generators of the $\Wc$-algebra
$\Wc^{\tss k}(\g,f)$. Moreover, we calculate the images of these generators with respect to
the {\it Miura transformation\/}. In particular, if $f$ is the principal nilpotent
(i.e., $n=1$) we thus reproduce the description of the $\Wc$-algebra due
to Fateev and Lukyanov~\cite{fl:mt}.
The results can be regarded as `affine analogues' of the construction of the corresponding
{\it finite $\Wc$-algebras\/} originated in \cite{br:rp}, \cite{rs:yr} and extended
to arbitrary nilpotent elements $f$ in \cite{bk:sy}.

\section{Principal $\Wc$-algebras}
\label{sec:pnp}
\setcounter{equation}{0}

The case where $f$ is a principal nilpotent will be useful in understanding
arbitrary rectangular $\Wc$-algebras so we consider it first.
We let $e_{ij}$ denote the standard basis elements of $\g=\gl_N$
and introduce subalgebras of $\g$ by
\ben
\finb
=\text{span of }\{e_{ij}\ |\ i\geqslant j\},\qquad
\finm
=\text{span of }\{e_{ij}\ |\ i> j\}\Fand
\finl=\text{span of }\{e_{ii}\}.
\een
Given $k\in\CC$ consider the normalized Killing form on $\g$ defined by
\beql{eq:form}
 \kappa(x,y)=\frac{k}{2 N}\on{tr}_{\fing}(\ad x\ts \ad y)
=k\ts\big(\tr(xy)-\frac{1}{N}\on{tr} x\ts \on{tr} y \big),\qquad x,y\in\g.
\eeq
Let
$
 \affb=\finb[t,t^{-1}]\+ \CC \mathbf{1}
$
be the Kac--Moody affinization of $\finb$ with respect to the
symmetric invariant bilinear form $\ka_{\tss\rm b}(x,y)$
on $\finb$, which is induced by \eqref{eq:form} and
given explicitly by
\ben
 \ka_{\tss\rm b}(e_{i\tss i\pr},\ts e_{j\tss j\pr})
=\delta_{i\tss  i\pr}\delta_{j\tss  j\pr}
\big(k+ N\big)\Big(\delta_{i\tss j}
-\frac{1}{N}\Big)
\een
for $i\geqslant i\pr$ and $j\geqslant j\pr$; see \cite[Lemma~4.8.1]{a:rtw}.
Furthermore,
let $V^{k}(\finb)$ be the universal affine vertex algebra
associated with $\finb$ and $\ka_{\tss\rm b}$ \cite{k:va}:
\beql{vkb}
 V^{k}(\finb)=\U(\affb)\ot_{\U(\finb[t]\+ \CC \mathbf{1})}\CC,
\eeq
where
$\CC$ is regarded as the one-dimensional
representation of
$\finb[t]\+ \CC \mathbf{1}$ on which $\finb[t]$ acts trivially
and $\mathbf{1}$ acts as $1$.
By the Poincar\'e--Birkhoff--Witt theorem,
$V^{k}(\finb)$ is isomorphic to $\U(\finb[t^{-1}]t^{-1})$ as a vector space.
Choose the principal nilpotent element in the form
\ben
f=\sum_{i=1}^{N-1}e_{i+1\ts i}\in \finm.
\een
The principal $\Wc$-algebra
$\Wc^{\tss k}(\fing)=\Wc^{\tss k}(\fing,f)$
can be realized as a vertex subalgebra of $V^{k}(\finb)$; see e.g.
\cite{a:rtw} and \cite[Ch.~15]{fb:va}.
Our aim is to give an explicit description of the generators
of
$\Wc^{\tss k}(\fing)$
inside $V^{k}(\finb)$.

Recall that the
{\em column-determinant} of a matrix $A=[a_{ij}]$ over
an associative algebra is defined by
\beql{cdet}
\on{cdet} A=\sum_{\sigma\in\Sym_N}\sgn\si\cdot a^{}_{\sigma(1)\tss 1}
a^{}_{\sigma(2)\tss 2}\dots\ts a^{}_{\sigma (N)\ts N}.
\eeq
Introduce an extended Lie algebra $\affb\oplus\CC\tau$, where the element $\tau$ commutes
with $\mathbf{1}$, and
\beql{tau}
\big[\tau, x[r]\big]=-r\tss x[r-1]\qquad\text{for}\quad
x\in \finb\fand
r\in \ZZ,
\eeq
where $x[r]=x\tss t^r$. In particular, \eqref{tau} induces
an associative algebra structure on
the tensor product
space $\U\big(\finb[t^{-1}]t^{-1}\big)\ot \CC[\tau]$.
Set $\alpha=k+N-1$
and consider the matrix
\begin{align*}
B=\begin{bmatrix}
\alpha\tss\tau+e_{11}[-1] &-1\phantom{-}&0&\dots & 0\\[0.4em]
e_{2\tss 1}[-1] &\alpha\tss\tau+e_{2\tss 2}[-1] &-1\phantom{-}&\dots & 0\\[0.4em]
\vdots &\vdots &\ddots & &\vdots\\[0.4em]
e_{N-1\tss 1}[-1] &e_{N-2\tss 2}[-1] &\dots
&\alpha\tss\tau+e_{N-1\ts N-1}[-1] &-1\phantom{-}\\[0.4em]
e_{N\tss 1}[-1] &e_{N\tss 2}[-1] &\dots &e_{N\tss N-1}[-1]  &\alpha\tss\tau+e_{N\tss N}[-1]
   \end{bmatrix}
\end{align*}
with entries in $\U\big(\finb[t^{-1}]t^{-1}\big)\ot \CC[\tau]$.
For its column-determinant\footnote{It is easy to verify that
$\on{cdet}B$ coincides with the {\it row-determinant\/}
of $B$ defined in a similar way.}
we can write
\begin{align*}
\on{cdet}B=(\alpha\tss\tau)^{\tss N}+W^{(1)}\ts(\alpha\tss\tau)^{\tss N-1}+\dots+
W^{(N)}
\end{align*}
for certain coefficients
$W^{(r)}$ which are elements of $\U\big(\finb[t^{-1}]t^{-1}\big)$,
and we can also regard them
as elements of $V^{k}(\finb)$.
We can now state our main result
on principal $\Wc$-algebras.

\bth\label{Th:mainp}
All coefficients $W^{(1)},\dots,W^{(N)}$ belong to $\Wc^{\tss k}(\fing)$.
Moreover, they freely generate the $\Wc$-algebra
$\Wc^{\tss k}(\fing)\subset V^{k}(\finb)$.
\eth

Before proving the theorem, note that
the projection
$\finb\ra \mf{l}$
induces the vertex algebra
homomorphism
$V^{k}(\finb)\ra V^{k}(\finl)$,
which restricts to the
map
\begin{align*}
 \nu:\Wc^{\tss k}(\fing)\ra   V^{k}(\mf{l}),
\end{align*}
called the ({\it quantum\/}) {\it Miura transformation}.
This is an injective
vertex algebra homomorphism; see \cite[Theorem~5.17]{a:iw}.
The following formula for the images of the elements $W^{(r)}$
under the Miura transformation is an immediate consequence
of Theorem~\ref{Th:mainp}. It reproduces
the construction of the $\Wc$-algebra $\Wc^{\tss k}(\fing)$
due to Fateev and Lukyanov~\cite{fl:mt}.

\bco\label{cor:fl}
Under the Miura transformation we have
\ben
\sum_{r=0}^{N}\nu(W^{(r)})(\alpha\tss \tau)^{N-r}
=\big(\alpha\tss\tau+e_{1\tss 1}[-1]\big)
\dots \big(\alpha\tss\tau+e_{N\tss N}[-1]\big).
\vspace{-0.7cm}
\een
\qed
\eco

\bpf[Proof of Theorem~\ref{Th:mainp}]
Taking into account \cite[Lemma~4.8.1]{a:rtw},
introduce the Lie superalgebra
$\affa=\affa_{\tss 0}\+ \affa_{1}$
such that  $\affa_{\tss 0}=\affb$ is the Lie subalgebra of even elements,
and $\affa_{1}=\finm[t,t^{-1}]$
is regarded as a supercommutative Lie superalgebra
spanned by odd elements,
while
\begin{align*}
 [x, y]=\ad x(y)\qquad\text{for}\quad x\in \affa_{\tss 0}\fand y\in \affa_{1}.
\end{align*}
We will use the notation
$\psi_{j\tss i}[m]=e_{j\tss i}\ts t^{m-1}$ (with a standard shift of the power by $1$)
for the element $e_{j\tss i}\ts t^{m-1}\in \finm[t,t^{-1}]$
when it is considered as an element of $\affa_{1}$.

Let
 $ V^{k}(\mf{a})$
be the representation of $\affa$ induced from the
one-dimensional representation of
$(\finb[t]
\+ \CC \mathbf{1})\+ \finm[t]$
on which
$\finb[t]\subset \affa_0$ and $\finm[t]\subset \affa_1$
act trivially and $\mathbf{1}$ acts as $1$.
Then $V^{k}(\mf{a})$ is naturally a vertex algebra
which contains $V^{k}(\finb)$ as its vertex subalgebra.
We will regard $V^{k}(\mf{a})$ as a (non-associative) algebra
with respect to the $(-1)$-product
\beql{mone}
 V^{k}(\mf{a})\ot V^{k}(\mf{a})\ra
 V^{k}(\mf{a}),
\qquad a\ot b\mapsto a_{(-1)}b,
\eeq
where the Fourier coefficients $a_{(n)}$ are defined in the usual way
from the state-field correspondence map,
\begin{align*}
 Y(a,z)=\sum_{n\in \ZZ}a_{(n)}z^{-n-1}\qquad\text{for}
 \quad a\in V^{k}(\mf{a}).
\end{align*}
Note
the quasi-associativity property of the $(-1)$-product
\beql{qass}
(a_{(-1)}b)_{(-1)}c=a_{(-1)}(b_{(-1)}c)
+\sum_{j\geqslant 0}a_{(-j-2)}(b_{(j)}c)+\sum_{j\geqslant 0}b_{(-j-2)}(a_{(j)}c)
\eeq
which holds for an arbitrary vertex algebra; see e.g.~\cite[Ch.~4]{k:va}.

The definition of the $\Wc$-algebra via the BRST cohomology \cite[Ch.~15]{fb:va}
can be stated in the form
\ben
 \Wc^{\tss k}(\fing)=\{v\in V^{k}(\mf{b})\ |\  Q\ts v=0\},
\een
where
$Q:V^{k}(\mf{a})\ra V^{k}(\mf{a})$ is a derivation
of the
non-associative algebra
$V^{k}(\mf{a})$
determined by the following properties. First, $Q$ commutes with the
translation operator $D$ of the vertex algebra $V^{k}(\mf{a})$.
Moreover, we have the commutation relations
\ben
\big[Q,e_{j\tss i}[-1]\big]=\sum_{a=i}^{j-1}
e_{a\tss i}[-1]\ts\psi_{j\tss a}[0]-\sum_{a=i+1}^j
\psi_{a\tss i}[0]\ts e_{j\tss a}[-1] +\alpha\ts \psi_{j\tss i}[-1]
+\psi_{j+1\ts i}[0] -\psi_{j\ts i-1}[0],
\een
which hold for all $j\geqslant i$ (empty sums are understood as being equal to zero),
and
\ben
\big[Q,\psi_{j\tss i}[0]\big]=\frac{1}{2}\sum_{i<r<j}
\big(\psi_{j\tss r}[0]\ts\psi_{r\tss i}[0]-\psi_{r\tss i}[0]\ts\psi_{j\tss r}[0]\big),
\een
for $j>i$, where $\psi_{k\tss l}[r]$ is considered to be zero for out-of-range indices,
and we omit the subscripts for the $(-1)$-products.

We will regard $V^{k}(\mf{a})\ot \CC[\tau]$ as a non-associative
algebra with the natural subalgebras $V^{k}(\mf{a})$ and
$\CC[\tau]$ together with the relation
$[\tau,u]=D\ts u$ for $u\in V^{k}(\mf{a})$.
In particular, the action of $Q$ on the extended algebra
commutes with the multiplication by $\tau$.

Now observe that the column-determinant $\cdet\ts B$ which we defined
over the associative algebra $\U\big(\finb[t^{-1}]t^{-1}\big)\ot \CC[\tau]$
coincides with the column-determinant $\wt{\on{cdet}}\ts B$
of the same matrix, but with the entries of $B$ regarded as elements of
the non-associative algebra $V^{k}(\mf{a})\ot \CC[\tau]$. Here we extend
the definition of
column-determinant to matrices $A=[a_{ij}]$ with entries
in a non-associative algebra by using right-normalized products,
\beql{nadet}
 \wt{\on{cdet}}\ts A=\sum_{\sigma\in\Sym_N}\sgn\si\cdot a^{}_{\sigma(1)\tss 1}
(a^{}_{\sigma(2)\tss 2}(a^{}_{\sigma (3)\tss 3}
\dots (a^{}_{\sigma(N-1)\ts  N-1}\ts a^{}_{\sigma (N)\ts N})\dots)).
\eeq
Thus, the first part of the theorem will follow from the relation
$
[Q, \wt{\cdet}\ts B]=0.
$
Expanding the column-determinant along the first column we get
\beql{cdetexpa}
\wt{\cdet}\ts B=\sum_{i=1}^N\big(\de_{i\tss 1}\tss\al\tss\tau+e_{i\tss 1}[-1]\big)\ts D^{(N-i)},
\eeq
where $D^{(p)}$ denotes the column-determinant of the $p\times p$
submatrix of $B$ corresponding to the last $p$ rows and columns (assuming
$D^{(0)}=1$). Using the induction on $N$, we derive
from the commutation relations satisfied by $Q$ that
\ben
[Q, D^{(N-i)}]=-\sum_{j=i+1}^N \psi_{j\tss i}[0]\ts D^{(N-j)}.
\een
Hence, by \eqref{cdetexpa} we can write
\begin{multline}
[Q, \wt{\cdet}\ts B]=\sum_{i=1}^{N}\Big(\sum_{a=1}^{i-1}e_{a\tss 1}[-1]\ts\psi_{i\tss a}[0]-\sum_{a=2}^i
 \psi_{a\tss 1}[0]\ts e_{i\tss a}[-1]
+\alpha\ts\psi_{i\tss 1}[-1]+\psi_{i+1\ts 1}[0]\Big)D^{(N-i)}\\
{}-\sum_{i=1}^{N}\big(\de_{i\tss 1}\tss\al\tss\tau+e_{i\tss 1}[-1]\big)
\ts\sum_{j=i+1}^N \psi_{j\tss i}[0]\ts D^{(N-j)}.
\non
\end{multline}
Apply the quasi-associativity \eqref{qass} in the first line
and change the order of summation in the second line to bring this
expression to the form
\ben
-\sum_{i=2}^{N}\sum_{a=2}^i
 \psi_{a\tss 1}[0]\ts e_{i\tss a}[-1]\ts D^{(N-i)}
+\sum_{i=2}^{N}\psi_{i\ts 1}[0]\ts D^{(N-i+1)}
-\al\ts\sum_{i=2}^{N} \psi_{i\tss 1}[0]\ts\tau\ts D^{(N-i)},
\een
where we used the relation
$\big[\tau, \psi_{j\tss i}[0]\big]=\psi_{j\tss i}[-1]$.
Therefore, changing the order of summation in the first term, we obtain
\ben
[Q, \wt{\cdet}\ts B]=\sum_{i=2}^{N}\psi_{i\ts 1}[0]\ts \Big(D^{(N-i+1)}
-\sum_{j=i}^{N}\big(\de_{j\tss i}\tss\al\tss\tau+e_{j\tss i}[-1]\big)
 \ts D^{(N-j)}\Big)=0,
\een
by applying the expansion along the first column for $D^{(N-i+1)}$.

The second part of Theorem~\ref{Th:mainp} follows by
considering the grading of $V^{k}(\finb)$ induced by the principal
grading of $\finb$ defined by $\deg e_{ij}=j-i$.
We have
\begin{align*}
W^{(r)}=\sum_{s=1}^{N-r+1}{e_{r+s-1\ts s}[-1]}
+\text{\ terms of higher degree}.
\end{align*}
The elements
$\sum_{s=1}^{N-r+1}{e_{r+s-1\ts s}}$ with $r=1,\dots,N$ form a basis of
the centralizer
$\fing^f$.
Hence the argument is completed by applying
\cite[Theorem 4.1]{kw:qr}; see also \cite[Theorem 5.5.1]{a:rt}.
\epf

\bre\label{rem:clwalg}
The {\em classical $\Wc$-algebra} $\Wc(\fing)$ can be obtained from $\Wc^{\tss k}(\fing)$
in three different ways: by taking the limit $k\to \infty$,
by taking the {\em critical level}
$k=-N$ \cite[Ch.~15]{fb:va},
or by considering
the Poisson vertex algebra associated with a canonical filtration~\cite[Sec.~5.10]{a:iw}.
The first way reproduces the classical
Miura transformation, whereas the second yields explicit generators
of $\Wc(\fing)$ inside the Poisson vertex algebra $V^{-N}(\finb)$; see also
\cite{mr:cw} for another proof and extensions to simple Lie algebras $\g$ of other types.
\qed
\ere

\section{Rectangular $\Wc$-algebras}
\label{sec:rwa}
\setcounter{equation}{0}

Now we consider the general rectangular case where
the Jordan type of a nilpotent element
$f\in\gl_N$ is of the form $(l^{\tss n})$ with $n\tss l=N$.
We identify $\g=\gl_N$ with the tensor product
of $\gl_{\tss l}$ and $\gl_n$ via the isomorphism $\gl_{\tss l}\ot\gl_n\to\g$ defined by
\beql{eq:N=nl}
e_{i\tss j}\ot e_{r\tss s}\mapsto  e_{(i-1)\tss n+r,\ts (j-1)\tss n+s},
\eeq
where the $e_{i\tss j}$ denote the standard basis elements of the
corresponding general linear Lie
algebras. Set
\ben
f_l=\sum_{i=1}^{l-1}e_{i+1\ts i}\in \gl_{\tss l}
\een
and
\ben
f=f_l\ot I_n=\sum_{i=1}^{l-1}\sum_{j=1}^n e_{i\tss n+j,\ts(i-1)\tss n+j}\in \g,
\een
where $I_n\in\gl_n$ is the identity matrix.
The matrix $f$ is a nilpotent element of $\g$
of Jordan type $(l^{\tss n})$.
Let
\ben
\gl_{\tss l}=\bigoplus_{p\in \ZZ}(\gl_{\tss l})_p
\een
be the standard
principal grading
of $\gl_{\tss l}$,
obtained by defining
the degree of $e_{i\tss j}$ to be equal to $j-i$.
Set
\ben
\gl_{\tss l,\leqslant 0}=\bigoplus_{p\leqslant 0}(\gl_{\tss l})_p\Fand
\gl_{\tss l,< 0}=\bigoplus_{p<0}(\gl_{\tss l})_p.
\een
The isomorphism \eqref{eq:N=nl}
then induces the $\ZZ$-grading on $\g$,
\ben
\g=\bigoplus_{p\in \ZZ}\g_p,\qquad
\g_p=(\gl_{\tss l})_p\ot \gl_n,
\een
which is a {\it good grading\/} for $f$ in the sense of
\cite{krw:qr}.
We also set
\beql{bm}
\finb
=\bigoplus_{p\leqslant 0}\g_p=\gl_{\tss l,\leqslant 0}\ot \gl_n\Fand
\finm
=\bigoplus_{p<0}\g_p=\gl_{\tss l,<0}\ot \gl_n.
\eeq

For any $k\in\CC$, consider the
symmetric invariant bilinear form on $\g$ defined in \eqref{eq:form}
and for elements $x,y\in\finb$ set
\begin{align*}
\ka_{\tss\rm b}(x,y)=
		     \kappa(x,y)+\frac{1}{2}\on{tr}_{\fing}(\ad x\ts \ad y)
-\frac{1}{2}\on{tr}_{\fing_0}p_0(\ad x\ts\ad y),
\end{align*}
where $p_0$ denotes the restriction of the operator to $\fing_0$.
Then $\ka_{\tss\rm b}$ defines
a symmetric invariant bilinear form
on $\finb$.
For $i\geqslant i\pr$ and $j\geqslant j\pr$ we have
\begin{multline}
 \ka_{\tss\rm b}(e_{i\tss i\pr}\ot e_{p\tss q},\ts e_{j\tss j\pr}\ot e_{r\tss s})\\
{}=\delta_{i\tss  i\pr}\delta_{j\tss  j\pr}
\Big((k+ n\tss l)\big(\delta_{i\tss j}\tss\delta_{p\tss s}\tss\delta_{q\tss r}
-\frac{1}{n\tss l}\tss\delta_{p\tss q}\tss\delta_{r\tss s}\big)
-n\tss \delta_{i\tss j}\big(\delta_{p\tss s}\tss\delta_{q\tss r}-\frac{1}{n}
\tss\delta_{p\tss q}\tss\delta_{r\tss s}\big)\Big).
\non
\end{multline}

As with the principal case, we let
$
 \affb=\finb[t,t^{-1}]\+ \CC \mathbf{1}
$
be the Kac--Moody affinization of $\finb$ with respect to the
form $\ka_{\tss\rm b}$,
and
let $V^{k}(\finb)$ be the universal affine vertex algebra
associated with $\finb$ and $\ka_{\tss\rm b}$ defined
as in \eqref{vkb}; see \cite{k:va}.
We have a vector space isomorphism
$V^{k}(\finb)\cong\U(\finb[t^{-1}]t^{-1})$.
Due to \cite{kw:qr} and \cite{kw:cq} (see also \cite{a:rt}),
the $\Wc$-algebra
$\Wc^{\tss k}(\fing,f)$
can be realized as a vertex subalgebra of $V^{k}(\finb)$.
We will give a description of the generators
of $\Wc^{\tss k}(\fing,f)$
inside $V^{k}(\finb)$.
We will use
the identification
\ben
 \mf{gl}_{\tss l,\leqslant 0}[t^{-1}]t^{-1}
\ot \mf{gl}_n\cong \finb[t^{-1}]t^{-1},
\een
defined by
\ben
e_{j\tss i}[-m]\ot e_{p\tss q}\mapsto
e_{(j-1)\tss n+p,\ts (i-1)\tss n+q}[-m],\qquad m\geqslant 1,
\een
for $1\leqslant i\leqslant j\leqslant l$ and $1\leqslant p,q\leqslant n$,
where we write $x[r]=x\ts t^{\tss r}$ for any $r\in\ZZ$.

By analogy with \cite[Sec.~12]{bk:sy}, consider the tensor algebra
$\Tr\big(\mf{gl}_{\tss l,\leqslant 0}[t^{-1}]t^{-1}\big)$ of the vector space
$\mf{gl}_{\tss l,\leqslant 0}[t^{-1}]t^{-1}$ and let $M_n$ denote the
matrix algebra with the basis formed by the matrix units $e_{i\tss j}$,
$1\leqslant i,j\leqslant n$.
Define the algebra homomorphism
\begin{align*}
 \Tc: \Tr\big(\mf{gl}_{\tss l,\leqslant 0}[t^{-1}]t^{-1}\big)\ra M_n\ot
 \U\big(\finb[t^{-1}]t^{-1}\big),
\qquad x\mapsto \Tc(x)=\sum_{i,j=1}^n e_{i\tss j}\ot \Tc_{i\tss j}(x)
\end{align*}
by setting
\begin{align*}
\Tc_{i\tss j}(x)=x\ot e_{j\tss i}\in \mf{gl}_{\tss l,\leqslant 0}[t^{-1}]t^{-1}
\ot \mf{gl}_n=\finb[t^{-1}]t^{-1}
\end{align*}
for $x\in \fing_{\tss l,\leqslant 0}[t^{-1}]t^{-1}$.
By definition, for any $x,y\in \Tr(\mf{gl}_{\tss l,\leqslant 0}[t^{-1}]t^{-1})$
we have
\begin{align*}
\Tc_{i\tss j}(x y)=\sum_{r=1}^n \Tc_{i\tss r}(x)\tss\Tc_{r\tss j}(y).
\end{align*}
Let us equip the tensor product space
$\Tr(\fing_{\tss l,\leqslant 0}[t^{-1}]t^{-1})
\ot \CC[\tau]$
with an associative algebra structure in such a way that
the natural embeddings
\ben
\Tr(\fing_{\tss l,\leqslant 0}[t^{-1}]t^{-1})\hookrightarrow
\Tr(\fing_{\tss l,\leqslant 0}[t^{-1}]t^{-1})
\ot \CC[\tau]\Fand
\CC[\tau]
\hookrightarrow
\Tr(\fing_{\tss l,\leqslant 0}[t^{-1}]t^{-1})
\ot \CC[\tau]
\een
are algebra homomorphisms
and the generator $\tau$ satisfies the relations
\begin{align*}
\big[\tau, x[-m]\big]=m\tss x[-m-1]\qquad\text{for}\quad
x\in \fing_{\tss l,\leqslant 0}\fand
m\in \ZZ.
\end{align*}
Furthermore, the tensor product space $\U(\finb[t^{-1}]t^{-1})\ot \CC[\tau]$ will also
be considered as an associative algebra in a similar way.
We will extend
$\Tc$ to the algebra
homomorphism
\begin{align*}
\Tc: \Tr(\fing_{\tss l,\leqslant 0}[t^{-1}]t^{-1})\ot \CC[\tau]
\ra M_n\ot \U(\finb[t^{-1}]t^{-1})\ot \CC[\tau]
\end{align*}
by setting $\Tc_{i\tss j}(u \tss S)=\Tc_{i\tss j}(u)\tss S$
for $u\in \Tr(\fing_{\tss l,\leqslant 0}[t^{-1}]t^{-1})$ and
any polynomial $S\in\CC[\tau]$.

Set $\alpha=k+n\tss (l-1)$
and consider the matrix
\begin{align*}
B=\begin{bmatrix}
\alpha\tss\tau+e_{11}[-1] &-1\phantom{-}&0&\dots & 0\\[0.4em]
e_{2\tss 1}[-1] &\alpha\tss\tau+e_{2\tss 2}[-1] &-1\phantom{-}&\dots & 0\\[0.4em]
\vdots &\vdots &\ddots & &\vdots\\[0.4em]
e_{l-1\tss 1}[-1] &e_{l-2\tss 2}[-1] &\dots
&\alpha\tss\tau+e_{l-1\ts l-1}[-1] &-1\phantom{-}\\[0.4em]
e_{l\tss 1}[-1] &e_{l\tss 2}[-1] &\dots &e_{l\ts l-1}[-1]  &\alpha\tss\tau+e_{l\tss l}[-1]
   \end{bmatrix}
\end{align*}
with entries in $\Tr(\mf{gl}_{\tss l,\leqslant 0}[t^{-1}]t^{-1})\ot \CC[ \tau]$.
Its column-determinant
$\on{cdet}B$ is defined by \eqref{cdet}.
So $\on{cdet}B$ is an element of
$\Tr(\mf{gl}_{\tss l,\leqslant0 }[t^{-1}]t^{-1})\ot \CC[\tau]$ and
we can write
\begin{align*}
\Tc_{i\tss j}(\on{cdet}B)=\sum_{r=0}^{l}W_{i\tss j}^{(r)}(\alpha\tss\tau)^{\tss l-r}
\end{align*}
for certain coefficients
$W_{i\tss j}^{(r)}$ which are elements of $\U(\finb[t^{-1}]t^{-1})$,
and we can also regard them
as elements of $V^{k}(\finb)$. The following is our main result
generalizing Theorem~\ref{Th:mainp}.

\bth\label{Th:main}
All coefficients $W_{i\tss j}^{(r)}$ belong to the $\Wc$-algebra $\Wc^{\tss k}(\fing,f)$.
Moreover, the $\Wc$-algebra
$\Wc^{\tss k}(\fing,f)\subset V^{k}(\finb)$
is freely generated by
the elements $W_{i\tss j}^{(r)}$ with
$1\leqslant i,j\leqslant n$ and $r=1,2,\dots, l$.
\eth

Set $\mf{l}=(\mf{gl}_{\tss l})_0\ot \mf{gl}_n\subset \mf{gl}_N$.
As with the principal case,
the projection
$\finb\ra \mf{l}$
induces the vertex algebra
homomorphism
$V^{k}(\finb)\ra V^{k}(\finl)$,
which restricts to the {\it Miura transformation}
\begin{align*}
 \nu:\Wc^{\tss k}(\fing,f)\ra   V^{k}(\mf{l}).
\end{align*}
This is an injective
vertex algebra homomorphism; see \cite[Sec.~5.9]{a:iw}.
The following formula for the images of the elements $W_{i\tss j}^{(r)}$
under the Miura transformation is an immediate consequence
of Theorem~\ref{Th:main}.

\bco\label{cor:miurarec}
We have
\begin{align*}
\sum_{r=0}^{l}\nu(W_{i\tss j}^{(r)})(\alpha\tss \tau)^{l-r}
=\Tc_{i\tss j}\Big(\big(\alpha\tss\tau+e_{1\tss 1}[-1]\big)
\dots \big(\alpha\tau+e_{l\tss l}[-1]\big)\Big).
\end{align*}
\eco

The principal $\Wc$-algebra corresponds to the case
$n=1$ (and $N=l$) so that Corollary~\ref{cor:miurarec}
generalizes the Fateev--Lukyanov formula; see Corollary~\ref{cor:fl}.
For a concrete example in the non-principal case see Example~\ref{ex:nonpri} below.

\bpf[Proof of Theorem~\ref{Th:main}]
Recall the notation \eqref{bm} and let
$\affa=\affa_{\tss 0}\+ \affa_{1}$ be the Lie superalgebra
such that elements of $\affa_{\tss 0}=\affb$ are even,
and elements of $\affa_{1}=\finm[t,t^{-1}]$ are odd,
where $\finm[t,t^{-1}]$ is regarded as the supercommutative Lie superalgebra,
while
\begin{align*}
 [x, y]=\ad x(y)\qquad\text{for}\quad x\in \affa_{\tss 0}\fand y\in \affa_{1}.
\end{align*}
We will write
$\psi_{j\tss i}[m]\ot e_{p\tss q}$
for the element
\ben
e_{j\tss i}\tss t^{m-1}\ot e_{p\tss q}\in \mf{gl}_{\tss l,<0}[t,t^{-1}]
\ot \mf{gl}_n=\finm[t,t^{-1}]
\een
when it is considered as an element of $\affa_{1}$.

Let
 $ V^{k}(\mf{a})$
be the representation of $\affa$ induced from the
one-dimensional representation of
$(\finb[t]
\+ \CC \mathbf{1})\+ \finm[t]$
on which
$\finb[t]\subset \affa_0$ and $\finm[t]\subset \affa_1$
act trivially and $\mathbf{1}$ acts as $1$.
Then $V^{k}(\mf{a})$ is naturally a vertex algebra
which contains $V^{k}(\finb)$ as its vertex subalgebra.
Endow $V^{k}(\mf{a})$ with the $(-1)$-product defined as in \eqref{mone}
so that it gets a structure of
a non-associative algebra.
By \cite{kw:qr} and \cite{kw:cq} the $\Wc$-algebra is given by
\begin{align*}
 \Wc^{\tss k}(\fing,f)=\{v\in V^{k}(\mf{b})\ |\  Q\ts v=0\},
\end{align*}
where
 $Q:V^{k}(\mf{a})\ra V^{k}(\mf{a})$ is the derivation
of the
non-associative algebra
$V^{k}(\mf{a})$
defined by the following properties. The map $Q$ commutes with the
translation operator $D$ of the vertex algebra $V^{k}(\mf{a})$
and we have the commutation relations
\ben
\bal[]
\big[Q,e_{j\tss i}[-1]\ot e_{p\tss q}\big]
{}&=\sum_{a=i}^{j-1}\sum_{r=1}^n
\big(e_{a\tss i}[-1]\ot e_{r\tss q}\big)
\big(\psi_{j\tss a}[0]\ot e_{p\tss r}\big)\\
{}&-\sum_{a=i+1}^j
\sum_{r=1}^n\big(\psi_{a\tss i}[0]\ot e_{r\tss q}\big)
\big(e_{j\tss a}[-1]\ot e_{p\tss r}\big)\\[0.4em]
{}&+\alpha\ts\psi_{j\tss i}[-1]\ot e_{p\tss q}
+ \psi_{j+1\ts i}[0]\ot e_{p\tss q}-\psi_{j \ts i-1}[0]\ot e_{p\tss q}
\eal
\een
and
\ben
\bal[]
\big[Q,\psi_{j\tss i}[0]\ot e_{p\tss q}\big]&=\frac{1}{2}\sum_{i<r<j,\ts
1\leqslant s\leqslant n} \big(\psi_{j\tss r}[0]\ot
e_{q\tss  s}\big)\big(\psi_{r\tss i}[0]\ot e_{s\tss p}\big)\\[0.4em]
{}&-\frac{1}{2}\sum_{i<r<j,\ts
1\leqslant s\leqslant n}\big(\psi_{r\tss i}[0]\ot e_{s\tss p}\big)
\big(\psi_{j\tss r}[0]\ot e_{q\tss s}\big)
\eal
\een
with $\al=k+n\tss (l-1)$.
Here we assumed that $\psi_{k\tss l}[r]=0$ for out-of-range subscripts and
used the fact that
\begin{align*}
 \on{tr}_{\finm}p_+\big(\ad( e_{j\tss i}\ot e_{p\tss q})
 \ts\ad( e_{i\tss j}\ot e_{q\tss p})\big)=
n\tss (l+j-i-1)
\end{align*}
for $1\leqslant i<j\leqslant l$ and
$1\leqslant p,q\leqslant n$, where
$p_+$ denotes the restriction of the operator to $\finm$.

Our goal now is to reduce the calculations to the principal nilpotent case.
To this end, when $n=1$
we will write $\bar{\mf{a}}$ and $\bar{\mf{b}}$
respectively, instead of $\mf{a}$ and $\mf{b}$,
and replace $k$ with $\bar k=k+(n-1)(l-1)$ in \eqref{eq:form}.
Consequently,
$V^{\bar k}(\bar{\mf{a}})$
will denote the vertex algebra
$V^{k}(\mf{a})$
with $n=1$ (and $k$ replaced by $\bar k$).
We let
$\overline Q$ denote
the operator $Q$ for $V^{\bar k}(\bar{\mf{a}})$.
The commutation relations for $\overline Q$ are given in the proof
of Theorem~\ref{Th:mainp}.

We will regard $V^{k}(\mf{a})\ot \CC[\tau]$ as a non-associative
algebra with the natural subalgebras $V^{k}(\mf{a})$ and
$\CC[\tau]$ together with the relation
$[\tau,u]=D\ts u$ for $u\in V^{k}(\mf{a})$.
Similarly, the tensor product
$V^{\bar k}(\bar{\mf{a}})\ot \CC[\tau]$ will be regarded
as a non-associative algebra with the relation
$[\tau,u]=\overline D\ts u$ for $u\in V^{\bar k}(\bar{\mf{a}})$,
where $\overline D$ denotes the
translation operator of the vertex algebra
$V^{\bar k}(\bar{\mf{a}})$.
Define the non-associative algebra homomorphism
\begin{align*}
 \wt\Tc:V^{\bar k}(\bar{\mf{a}})\ot \CC[\tau]
 \ra M_n\ot V^{k}(\mf{a})\ot \CC[\tau],
\qquad x\mapsto \wt\Tc(x)=\sum_{p,q=1}^n e_{p\tss q}\ot \wt\Tc_{p\tss q}(x)
\end{align*}
by
\begin{align*}
 \wt\Tc_{p\tss q}\big(e_{j\tss i}[m]\big)=e_{j\tss i}[m]\ot e_{q\tss p},\qquad
 \wt\Tc_{p\tss q}\big(\psi_{j\tss i}[m]\big)=\psi_{j\tss i}[m]\ot e_{q\tss p}\Fand
\wt\Tc_{p\tss q}(\tau)=\tau.
\end{align*}

Note the relation
$
\wt\Tc(\wt{\on{cdet}}\ts B)=\Tc(\on{cdet}B),
$
where
$\wt{\on{cdet}}\ts B$ is defined by \eqref{nadet} and regarded as an element of
$V^{\bar k}(\bar{\mf{b}})\ot \CC[\tau]$.
The first part of
Theorem~\ref{Th:main} will now be implied by
the property
\ben
\big[Q,\wt\Tc_{p\tss q}(a)\big]=\wt\Tc_{p\tss q}\big([\overline Q,a]\big)
\een
which holds for any $a\in V^{\bar k}(\bar{\mf{b}})\ot \CC[\tau]$ and
$1\leqslant p,q\leqslant n$,
and which follows from
the definitions of the operators $Q$ and $\overline Q$.
It remains to note that the relation
$
[\overline Q, \wt{\on{cdet}}\ts B]=0
$
was already verified
in the proof of Theorem~\ref{Th:mainp}.

To prove the second part of Theorem~\ref{Th:main},
consider the grading of $V^{k}(\finb)$ induced by the
grading
of $\finb$.
One has
\begin{align*}
W_{ij}^{(r)}=\Tc_{ij}\Big(\sum_{s=1}^{l-r+1}{e_{r+s-1\ts s}[-1]}\Big)
+\text{terms of higher degree}.
\end{align*}
The elements
$\sum_{s=1}^{l-r+1}{e_{r+s-1\ts s}}$ with $r=1,\dots,l$ form a basis of
$\mf{gl}_{\tss l}^{f_l}$
and the elements
\ben
\sum_{s=1}^{l-r+1}{e_{r+s-1\ts s}}\ot e_{j\tss i},\qquad r=1,\dots,l\fand  i,j=1,\dots,n,
\een
form a basis of $\fing^f$.
It remains to apply
\cite[Theorem 4.1]{kw:qr}; see also \cite[Theorem 5.5.1]{a:rt}.
\epf

\bex\label{ex:nonpri}
Take $n=l=2$ so that $N=4$. We have
\begin{align*}
 \on{cdet}B=(\alpha\tss\tau)^2+\big(e_{1\tss 1}[-1]+e_{2\tss 2}[-1]\big)(\alpha\tss\tau)
 +e_{1\tss 1}[-1]\tss e_{2\tss 2}[-1]
+e_{2\tss 1}[-1]+\alpha\tss e_{2\tss 2}[-2]
\end{align*}
with $\al=k+2$. Hence
\begin{align*}
 &W_{1\tss 1}^{(1)}=e_{1\tss 1}[-1]+e_{3\tss 3}[-1],
\qquad  W_{2\tss 2}^{(1)}=e_{2\tss 2}[-1]+e_{4\tss 4}[-1],\\[0.3em]
& W_{2\tss 1}^{(1)}=e_{1\tss 2}[-1]+e_{3\tss 4}[-1],\qquad
 W_{1\tss 2}^{(1)}=e_{2\tss 1}[-1]+e_{4\tss 3}[-1],\\[0.3em]
&W_{1\tss 1}^{(2)}=e_{1\tss 1}[-1]\tss e_{3\tss 3}[-1]
+e_{2\tss 1}[-1]\tss e_{3\tss 4}[-1]+e_{3\tss 1}[-1]
+\alpha\ts e_{3\tss 3}[-2],\\[0.3em]
&W_{2\tss 2}^{(2)}=e_{1\tss 2}[-1]\tss e_{4\tss 3}[-1]
+e_{2\tss 2}[-1]\tss e_{4\tss 4}[-1]
+e_{4\tss 2}[-1]+\alpha\ts e_{4\tss 4}[-2],\\[0.3em]
 &W_{2\tss 1}^{(2)}=e_{1\tss 2}[-1]\tss e_{3\tss 3}[-1]
 +e_{2\tss 2}[-1]\tss e_{3\tss 4}[-1]
 +e_{3\tss 2}[-1]+\alpha\ts e_{3\tss 4}[-2],\\[0.3em]
 &W_{1\tss 2}^{(2)}=e_{1\tss 1}[-1]\tss e_{4\tss 3}[-1]
 +e_{2\tss 1}[-1]\tss e_{4\tss 4}[-1]
 +e_{4\tss 1}[-1]+\alpha\ts e_{4\tss 3}[-2].
\end{align*}

For the images under the Miura transformation we have
\begin{align*}
 &\nu(W_{1\tss 1}^{(1)})=e_{1\tss 1}[-1]+e_{3\tss 3}[-1],
\qquad  \nu(W_{2\tss 2}^{(1)})=e_{2\tss 2}[-1]+e_{4\tss 4}[-1],\\[0.3em]
& \nu(W_{2\tss 1}^{(1)})=e_{1\tss 2}[-1]+e_{3\tss 4}[-1],\qquad
 \nu(W_{1\tss 2}^{(1)})=e_{2\tss 1}[-1]+e_{4\tss 3}[-1],\\[0.3em]
&\nu(W_{1\tss 1}^{(2)})=e_{1\tss 1}[-1]\tss e_{3\tss 3}[-1]
+e_{2\tss 1}[-1]\tss e_{3\tss 4}[-1]
+\alpha\ts e_{3\tss 3}[-2],\\[0.3em]
&\nu(W_{2\tss 2}^{(2)})=e_{1\tss 2}[-1]\tss e_{4\tss 3}[-1]
+e_{2\tss 2}[-1]\tss e_{4\tss 4}[-1]
+\alpha\ts e_{4\tss 4}[-2],\\[0.3em]
 &\nu(W_{2\tss 1}^{(2)})=e_{1\tss 2}[-1]\tss e_{3\tss 3}[-1]
 +e_{2\tss 2}[-1]\tss e_{3\tss 4}[-1]
 +\alpha\ts e_{3\tss 4}[-2],\\[0.3em]
 &\nu(W_{1\tss 2}^{(2)})=e_{1\tss 1}[-1]\tss e_{4\tss 3}[-1]
 +e_{2\tss 1}[-1]\tss e_{4\tss 4}[-1]
 +\alpha\ts e_{4\tss 3}[-2].
\end{align*}

The values of the form $\ka_{\tss\rm b}(x,y)$ are given in the following table,
where the columns and rows correspond to the $x$ and $y$ variables,
respectively:

\bigskip

\begin{tabular}{|c|c|c|c|c|c|c|c|c|}
\hline
 &$e_{11}$ &$e_{22}$&$e_{33}$  &$e_{44}$ &$e_{12}$ &$e_{21}$ &$e_{34}$
 &$e_{43}$ \\
\hline
$e_{11}$ & $\frac{3k+8}{4}
 $ & $-\frac{k}{4}$& $-\frac{k+4}{4}
$ &$-\frac{k+4}{4}$ & $0$& $0$& $0$& $0$\\[0.2em]
\hline
$e_{22}$ & $-\frac{k}{4} $ & $\frac{3k+8}{4}
$& $-\frac{k+4}{4}$ &$-\frac{k+4}{4}$ & $0$& $0$& $0$& $0$\\[0.2em]
\hline
$e_{33}$ & $-\frac{k+4}{4} $ & $-\frac{k+4}{4}$& $\frac{3k+8}{4}$ &
$-\frac{k}{4}$ & $0$& $0$& $0$& $0$\\[0.2em]
\hline
$e_{44}$ & $-\frac{k+4}{4} $ & $-\frac{k+4}{4}$& $-\frac{k}{4}$ &
$\frac{3k+8}{4}$ & $0$& $0$& $0$& $0$
\\[0.2em]
\hline
$e_{12}$ & $0 $ & $0$& $0$ &$0$ & $0$& $k+2$& $0$& $0$
\\\hline
$e_{21}$ & $0 $ & $0$& $0$ &$0$ & $k+2$& $0$& $0$& $0$
\\\hline
$e_{34}$ & $0 $ & $0$& $0$ &$0$ & $0$& $0$& $0$& $k+2$
\\\hline
$e_{34}$ & $0 $ & $0$& $0$ &$0$ & $0$& $0$& $k+2$& $0$\\
\hline
\end{tabular}

\vspace{7mm}

\noindent
These values can be used to calculate the operator product expansion
formulas for the generators of $\Wc^{\tss k}(\fing,f)$. In particular,
set
\begin{align*}
 L=&\frac{1}{2(k+4)}
\Big({-2}(W_{11}^{(2)}+W_{22}^{(2)})
+W_{12}^{(1)}W_{21}^{(1)}
+\frac{3}{4}(W_{11}^{(1)}W_{11}^{(1)}+
W_{22}^{(1)}W_{22}^{(1)}
)\\[0.5em]
& -\frac{1}{2}W_{11}^{(1)}W_{22}^{(1)}
-(k+2)(W_{11}^{(1)}+W_{22}^{(1)})'-(W_{11}^{(1)}-W_{22}^{(1)})'\Big),
\end{align*}
where the primes indicate the action of $\ad\tau$ taking $e_{i\tss j}[-1]$
to $e_{i\tss j}[-2]$.
Then
$L$ is the conformal vector of $\Wc^{\tss k}(\fing,f)$:
\begin{align*}
 L(z)L(w)\sim -\frac{
12k^2+41 k+32}{2(k+4)^2(z-w)^4}
+\frac{2}{(z-w)^2}L(w)+\frac{1}{z-w}\partial L(w).
\end{align*}
\eex


\begin{thebibliography}{99}


\bibitem{a:rt}
T. Arakawa,
{\it Representation theory of superconformal algebras and
the Kac--Roan--Wakimoto conjecture},
 Duke Math. J. {\bf 130} (2005), 435--478.

\bibitem{a:rtw}
T. Arakawa,
{\it Representation theory of W-algebras},
Invent. Math. {\bf 169} (2007), 219--320.

\bibitem{a:rw}
T. Arakawa,
{\it Rationality of W-algebras: principal nilpotent cases},
Ann. of Math. (2) {\bf 182} (2015), 565--604.

\bibitem{a:iw}
T. Arakawa,
{\it Introduction to W-algebras and their representation theory},\\
{\tt arXiv:1605.00138}.

\bibitem{br:rp}
{C. Briot and E. Ragoucy},
{\it RTT presentation of finite W-algebras},
{J. Phys. A} {\bf 34} (2001), 7287--7310.

\bibitem{bk:sy}
{J. Brundan and A. Kleshchev},
{\it Shifted Yangians and finite W-algebras},
Adv. Math.  {\bf 200}  (2006),  136--195.

\bibitem{fl:mt}
V. A. Fateev and S. L. Lukyanov,
{\it The models of two-dimensional
conformal quantum field theory with $Z_n$ symmetry},
Internat. J. Modern Phys. A {\bf 3} (1988), 507--520.

\bibitem{ff:qd}
B. Feigin and E. Frenkel,
{\it Quantization of the Drinfeld--Sokolov reduction},
Phys. Lett. B {\bf 246} (1990), 75--81.

\bibitem{fb:va}
E. Frenkel and D. Ben-Zvi,
{\it Vertex algebras and algebraic curves},
Second edition.
Mathematical Surveys and Monographs, 88.
AMS, Providence, RI, 2004.

\bibitem{k:va}
V. Kac,
{\it Vertex algebras for beginners},
University Lecture Series, 10. American Mathematical Society,
Providence, RI, 1997.

\bibitem{krw:qr}
V. Kac, Shi-Shyr Roan and M. Wakimoto,
{\it Quantum reduction for affine superalgebras},
Comm. Math. Phys. {\bf 241} (2003), 307--342.

\bibitem{kw:qr}
V. Kac and M. Wakimoto,
{\it Quantum reduction
and representation theory of superconformal
algebras}, Adv. Math. {\bf 185} (2004), 400--458.

\bibitem{kw:cq}
V. Kac and M. Wakimoto,
{\it Corrigendum to:
``Quantum reduction and representation theory of
superconformal algebras" [Adv. Math. 185 (2004), 400--458]},
Adv. Math. {\bf 193} (2005), 453--455.

\bibitem{mr:cw}
A. I. Molev and E. Ragoucy,
{\it Classical W-algebras in types $A, B, C, D$ and $G$},
Comm. Math. Phys. {\bf 336} (2015), 1053--1084.

\bibitem{rs:yr}
{E. Ragoucy and P. Sorba},
{\it Yangian realisations from finite W-algebras},
{Comm. Math. Phys.} {\bf 203}
(1999), 551--572.



\end{thebibliography}
\end{document}